\documentclass[a4paper]{amsart}
\usepackage{color}
\usepackage[latin1]{inputenc}
\usepackage[T1]{fontenc}
\usepackage{amsfonts}
\usepackage{amssymb}
\usepackage{amsmath}
\usepackage{amsthm}
\usepackage{stmaryrd}
\usepackage{enumerate}
\usepackage{multirow}
\usepackage{graphicx}
\usepackage{comment}

\usepackage{graphicx,type1cm,eso-pic,color}
\usepackage{pstricks}
\usepackage{float}
\newtheorem{theorem}{Theorem}[section]
\newtheorem{lemma}[theorem]{Lemma}
\newtheorem{proposition}[theorem]{Proposition}

\newtheorem{assumption}[theorem]{Assumption}
\newtheorem{remark}[theorem]{Remark}

\begin{document}
\setlength\arraycolsep{2pt}
\title[Warped Kernel Estimator for I.I.D. Paths of Diffusion Processes]{Warped Kernel Estimator for I.I.D. Paths of Diffusion Processes}
\author{Nicolas MARIE$^{\dag}$}
\email{nmarie@parisnanterre.fr}
\author{Am\'elie ROSIER$^{\diamond,\ast}$}
\email{amelie.rosier@esme.fr}
\address{$^{\dag}$Laboratoire Modal'X, Universit\'e Paris Nanterre, Nanterre, France}
\address{$^{\diamond}$ESME, Ivry-sur-Seine, France}
\address{$^{\ast}$UMR MIA-Paris-Saclay, INRAE, France}
\date{}
\maketitle
%


%
\begin{abstract}
This paper deals with a nonparametric warped kernel estimator $\widehat b$ of the drift function computed from independent continuous observations of a diffusion process. A risk bound on $\widehat b$ is established. The paper also deals with an extension of the PCO bandwidth selection method for $\widehat b$. Finally, some numerical experiments are provided.
\end{abstract}
\textbf{MSC2020 subject classifications.} 62G05 ; 62M05.
\tableofcontents
%


%
\section{Introduction}\label{section_introduction}\label{section_introduction}
Consider $T > 0$ and the stochastic differential equation
\begin{equation}\label{main_equation}
X_t = x_0 +\int_{0}^{t}b(X_s)ds +\int_{0}^{t}\sigma(X_s)dW_s
\textrm{ $;$ }t\in [0,T],
\end{equation}
where $b,\sigma :\mathbb R\rightarrow\mathbb R$ are two continuous functions and $W = (W_t)_{t\in [0,T]}$ is a Brownian motion.
\\
\\
Since the 1980's, the statistical inference for stochastic differential equations (SDE) has been widely investigated by many authors in the parametric and in the nonparametric frameworks. Classically (see Hoffman \cite{HOFFMANN99}, Kessler \cite{KESSLER00}, Kutoyants \cite{KUTOYANTS04}, Dalalyan \cite{DALALYAN05}, Comte et al. \cite{CGCR07}, etc.), the estimators of the drift function are computed from one path of the stationary solution of Equation (\ref{main_equation}), which exists and is unique under a restrictive dissipativity condition on $b$, and converge when $T$ goes to infinity.
\\
\\
Let $\mathcal I(.)$ be the It\^o map for Equation (\ref{main_equation}) and, for $N\in\mathbb N^*$ copies $W^1,\dots,W^N$ of $W$, consider
\begin{displaymath}
X^i :=\mathcal I(x_0,W^i)
\textrm{ $;$ }
\forall i\in\{1,\dots,N\}.
\end{displaymath}
The estimation of the drift function $b$ from continuous-time and discrete-time observations of $(X^1,\dots,X^N)$ is a functional data analysis problem already investigated in the parametric framework (see Ditlevsen and De Gaetano \cite{DDG05}, Overgaard et al. \cite{OJTM05}, Picchini, De Gaetano and Ditlevsen \cite{PDGD10}, Picchini and Ditlevsen \cite{PD11}, Comte, Genon-Catalot and Samson \cite{CGCS13}, Delattre and Lavielle \cite{DL13}, Delattre, Genon-Catalot and Samson \cite{DGCS13}, Dion and Genon-Catalot \cite{DGC16}, Delattre, Genon-Catalot and Lar\'edo \cite{DGCL13}, etc.) and more recently in the nonparametric one (see Comte and Genon-Catalot \cite{CGC20b,CGC21}, Comte and Marie \cite{CM23}, Della Maestra and Hoffmann \cite{DMH21}, Marie and Rosier \cite{MR23}, Denis et al. \cite{DDM21}).
\\
\\
Under the appropriate conditions on $b$ and $\sigma$ recalled in Section \ref{section_risk_bound}, for every $t\in (0,T]$, the distribution of $X_t$ has a density $p_t(x_0,.)$ such that $s\mapsto p_s(x_0,x)$ belongs to $\mathbb L^1([0,T])$ for every $x\in\mathbb R$. Then, for a given $t_0\in [0,T]$, one may define
\begin{displaymath}
f(x) :=\frac{1}{T - t_0}\int_{t_0}^{T}p_t(x_0,x)dt
\quad {\rm and}\quad
F(x) :=\int_{-\infty}^{x}f(z)dz
\end{displaymath}
for every $x\in\mathbb R$. Clearly, $f$ is a density function:
\begin{displaymath}
\int_{-\infty}^{\infty}f(x)dx =
\frac{1}{T - t_0}\int_{t_0}^{T}\int_{-\infty}^{\infty}p_t(x_0,x)dxdt = 1.
\end{displaymath}
Let $K :\mathbb R\rightarrow\mathbb R$ be a kernel (i.e. an integrable function such that $\int K = 1$), and consider $K_h(x) := h^{-1}K(h^{-1}x)$ with $h\in (0,1]$. Our paper deals with the warped kernel estimator
\begin{displaymath}
\widehat b_{N,h}(x) :=\widehat\beta_{N,h}(F;F(x))
\textrm{ $;$ }
x\in\mathbb R
\end{displaymath}
of the drift function $b$, where
\begin{displaymath}
\widehat\beta_{N,h}(\varphi;z) :=
\frac{1}{N(T - t_0)}\sum_{i = 1}^{N}
\int_{t_0}^{T}K_h(z -\varphi(X_{t}^{i}))dX_{t}^{i}
\end{displaymath}
for every $\varphi\in C^0(\mathbb R)$ and $z\in\mathbb R$. From independent copies of $X$ continuously observed on $[0,T]$, $\widehat b_{N,h}$ is a natural extension of the warped kernel estimator already well-studied in the nonparametric regression framework (see Chagny \cite{CHAGNY15}). The paper also deals with the adaptive estimator
\begin{displaymath}
\widehat b_N(x) :=
\widehat b_{N,\widehat h}(x),
\end{displaymath}
where $\widehat h$ is selected via a penalized comparison to overfitting (PCO) type criterion. Finally, in practice, the function $F$ is unknown and has to be replaced by the estimator
\begin{displaymath}
\widehat F_N(x) :=
\frac{1}{N(T - t_0)}\sum_{i = 1}^{N}\int_{t_0}^{T}\mathbf 1_{X_{t}^{i}\leqslant x}dt
\end{displaymath}
in the definition of $\widehat b_{N,h}$. Precisely, $\widehat b_{N,h}$ is approximated by
\begin{displaymath}
\widetilde b_{N,h}(x) :=\widehat\beta_{N,h}(\widehat F_N;\widehat F_N(x)).
\end{displaymath}
Section \ref{section_risk_bound} deals with a risk bound on the warped kernel estimator and Section \ref{section_PCO} with a risk bound on the adaptive estimator $\widehat b_N$. Section \ref{section_numerical_experiments} deals with some numerical experiments on $\widetilde b_{N,h}$. The proofs are postponed to Section \ref{section_proofs}.
\\
\\
\textbf{Notations and basic definitions:}
\begin{itemize}
 \item The space $C^0(\mathbb R)$ is equipped with the uniform (semi-)norm $\|.\|_{\infty}$.
 \item For every $p\in\overline{\mathbb N}$, $C_{\rm b}^{p}(\mathbb R) :=\cap_{j = 0}^{p}\{\varphi\in C^p(\mathbb R) :\varphi^{(j)}\textrm{ is bounded}\}$.
 \item For every $p\geqslant 1$, $\mathbb L^p(\mathbb R,dx)$ is equipped with its usual norm $\|.\|_p$:
 \begin{displaymath}
 \|\varphi\|_p :=\left(\int_{-\infty}^{\infty}\varphi(x)^pdx\right)^{1/p}
 \textrm{$;$ }
 \forall\varphi\in\mathbb L^p(\mathbb R,dx).
 \end{displaymath}
 \item $\mathbb H^2$ is the space of the processes $(Y_t)_{t\in [0,T]}$, adapted to the filtration generated by $W$, such that
 \begin{displaymath}
 \int_{0}^{T}\mathbb E(Y_{t}^{2})dt <\infty.
 \end{displaymath}
 \item For a given kernel $\delta$, the usual scalar product on $\mathbb L^2(\mathbb R,\delta(x)dx)$ is denoted by $\langle .,.\rangle_{\delta}$, and the associated norm by $\|.\|_{\delta}$. For every $A,B\in\overline{\mathbb R}$ such that $A < B$, the usual norm on $\mathbb L^2([A,B],\delta(x)dx)$ is denoted by $\|.\|_{\delta,A,B}$.
\end{itemize}
%


%
\section{Risk bound on the warped kernel estimator}\label{section_risk_bound}
In the sequel, in order to ensure the existence and the uniqueness of the (strong) solution of Equation (\ref{main_equation}), and to ensure that the distribution of $X_t$ has a regular enough density $p_t(x_0,.)$ for every $t\in (0,T]$ (see Remark \ref{F_properties}), $b$ and $\sigma$ fulfill the following assumption.
%


%
\begin{assumption}\label{assumption_b_sigma}
The function $b$ is Lipschitz continuous, $\sigma\in C_{\rm b}^{1}(\mathbb R)$, $\sigma'$ is H\"older continuous and there exists $\alpha > 0$ such that
\begin{displaymath}
|\sigma(x)| >\alpha
\textrm{ $;$ }
\forall x\in\mathbb R.
\end{displaymath}
\end{assumption}
%


%
\begin{remark}\label{F_properties}
By Menozzi et al. \cite{MPZ21}, Theorem 1.2, for any $t\in (0,T]$, the distribution of $X_t$ has a continuously differentiable and sub-Gaussian density $p_t(x_0,.)$ such that $p_t(x_0,\mathbb R)\subset (0,\infty)$, $z\in\mathbb R\mapsto\partial_zp_t(x_0,z)$ is also sub-Gaussian and $s\in (0,T]\mapsto p_s(x_0,z)$ belongs to $\mathbb L^1([0,T],dt)$ for every $z\in\mathbb R$. First, by this last point,
\begin{displaymath}
F(.) =\int_{-\infty}^{.}f(z)dz =
\frac{1}{T - t_0}\int_{t_0}^{T}F_s(.)ds
\quad {\rm with}\quad
F_s(.) :=
\int_{-\infty}^{.}p_s(x_0,z)dz
\textrm{ $;$ }
\forall s\in (0,T]
\end{displaymath}
is well-defined (even when $t_0 = 0$). Moreover, $F_t$ (and then $F$) is strictly increasing and one-to-one from $\mathbb R$ into $(0,1)$ because $p_t(x_0,.)$ is continuous and $p_t(x_0,\mathbb R)\subset (0,\infty)$. Finally, since $p_t(x_0,.)$ is sub-Gaussian, $b$ and $\sigma$ belong to $\mathbb L^2(\mathbb R,f(z)dz)$.
\end{remark}
\noindent
The kernel $K$ fulfills the following usual assumption.
%


%
\begin{assumption}\label{assumption_K}
The kernel $K$ is symmetric, continuous and belongs to $\mathbb L^2(\mathbb R,dx)$.
\end{assumption}
\noindent
The following proposition deals with a risk bound on $\widehat b_{N,h}$.
%


%
\begin{proposition}\label{risk_bound_warped_estimator}
Consider $A,B\in\overline{\mathbb R}$ such that $A < B$. Under Assumptions \ref{assumption_b_sigma} and \ref{assumption_K},
\begin{displaymath}
\mathbb E(\|\widehat b_{N,h} - b\|_{f,A,B}^{2})
\leqslant
\|b_h - b\|_{f,A,B}^{2} +
\frac{\mathfrak c_{\ref{risk_bound_warped_estimator}}}{Nh}
\end{displaymath}
with
\begin{displaymath}
b_h :=
[K_h\ast ((b\circ F^{-1})\mathbf 1_{(0,1)})](F(.))
\quad {\rm and}\quad
\mathfrak c_{\ref{risk_bound_warped_estimator}} =
2\|K\|_{2}^{2}
\left(\|b\|_{f}^{2} + \frac{1}{T - t_0}\|\sigma\|_{f}^{2}\right).
\end{displaymath}
\end{proposition}
\noindent
Let us conclude this section with some remarks about Proposition \ref{risk_bound_warped_estimator}:
\begin{itemize}
 \item Note that the variance term in the risk bound on $\widehat b_{N,h}$ stated in Proposition \ref{risk_bound_warped_estimator} is of same order than in the nonparametric regression framework (see Chagny \cite{CHAGNY15}).
 \item If $b\in\mathbb L^2(\mathbb R,dx)$, then
 \begin{displaymath}
 b_h\xrightarrow[h\rightarrow 0^+]{\mathbb L^2}
 (b\circ F^{-1})(F(.)) = b.
 \end{displaymath}
 \item Assume that $K$ is a $[-1,1]$-supported kernel, $A,B\in\mathbb R$ and that $h\in (0,\kappa(F(A)\wedge (1 - F(B)))]$ with $\kappa\in (0,1)$. Then, for every $x\in [A,B]$,
 \begin{displaymath}
 [-1,1]\subset
 \left[-\frac{F(A)}{h},\frac{1 - F(B)}{h}\right]\subset
 \left[-\frac{F(x)}{h},\frac{1 - F(x)}{h}\right]
 \end{displaymath}
 and
 \begin{displaymath}
 F^{-1}(h\cdot + F(x))([-1,1])
 \subset
 [F^{-1}((1 -\kappa)F(A)),F^{-1}(\kappa + (1 -\kappa)F(B))] =: I_{A,B,F}.
 \end{displaymath}
 Moreover, since $f(\mathbb R)\subset (0,\infty)$, there exists $f_1 > 0$ such that $f(z)\geqslant f_1$ for every $z\in I_{A,B,F}$. So,
 \begin{eqnarray*}
  \|b_h - b\|_{f,A,B}^{2} & = &
  \int_{A}^{B}\left(
  \int_{0}^{1}K_h(y - F(x))(b\circ F^{-1})(y)dy - b(x)\right)^2f(x)dx\\
  & = &
  \int_{A}^{B}\left[\int_{-F(x)/h}^{(1 - F(x))/h}
  K(y)b(F^{-1}(hy + F(x)))dy - b(x)\right]^2f(x)dx\\
  & = &
  \int_{A}^{B}\left(\int_{-1}^{1}K(y)
  [b(F^{-1}(hy + F(x))) - b(F^{-1}(F(x)))]dy\right)^2f(x)dx\\
  & \leqslant &
  h^2\|b\|_{\rm Lip}^{2}
  \left(\sup_{z\in I_{A,B,F}}\frac{1}{(f\circ F^{-1})(z)}\right)^2\\
  & &
  \hspace{1.5cm}\times
  \int_{-1}^{1}K(y)y^2\int_{A}^{B}f(x)dxdy\leqslant
  \mathfrak c_1h^2
  \quad {\rm with}\quad
  \mathfrak c_1 =
  \frac{\|b\|_{\rm Lip}^{2}}{f_{1}^{2}}\int_{-1}^{1}K(y)y^2dy.
 \end{eqnarray*}
 Therefore, the bias term in the risk bound on $\widehat b_{N,h}$ stated in Proposition \ref{risk_bound_warped_estimator} is of order $h^2$, as in the nonparametric regression framework (see Chagny \cite{CHAGNY15}).
\end{itemize}
%


%
\section{A PCO type bandwidth selection method}\label{section_PCO}
This section deals with an extension of the PCO method to the warped kernel estimator studied in this paper. The PCO based adaptive warped kernel estimator offers theoretical guarantees: a risk bound is established in this section.
\\
\\
In this section, consider an additional kernel $\delta$ and assume that $K$, $\delta$ and $\sigma$ fulfill the following technical assumption.
%


%
\begin{assumption}\label{assumption_K_sigma_PCO}
The kernels $K$ and $\delta$ are compactly supported and continuously differentiable on $\mathbb R$, $\delta(I)\subset (0,\infty)$ for every compact subset $I$ of ${\rm supp}(\delta)^{\circ}$, and $\sigma$ is bounded.
\end{assumption}
\noindent
For instance, the kernel
\begin{displaymath}
\rho : x\in\mathbb R\longmapsto
\frac{1}{\mathfrak c_{\rho}}\exp\left(-\frac{1}{1 - x^2}\right)\mathbf 1_{[-1,1]}(x)
\quad {\rm with}\quad
\mathfrak c_{\rho} =
\int_{-1}^{1}\exp\left(-\frac{1}{1 - y^2}\right)dy
\end{displaymath}
satisfies the conditions on both $K$ and $\delta$ in Assumption \ref{assumption_K_sigma_PCO}. In the sequel, for the sake of simplicity, ${\rm supp}(K) = [-1,1]$ and ${\rm supp}(\delta) = [-\Delta,\Delta]$ with $\Delta > 0$. Now, let $\mathcal H_N$ be a finite subset of $[h_0,\Delta_0]$, where
\begin{displaymath}
\Delta_0 :=
\kappa
(F(-\Delta)\wedge (1 - F(\Delta))),
\quad\kappa\in (0,1),
\quad h_0\in (0,\Delta_0)
\quad {\rm and}\quad
\frac{\Delta_{0}^{3}}{Nh_{0}^{3}}\leqslant 1.
\end{displaymath}
Finally, consider
\begin{equation}\label{penalty_proposal}
\widehat h\in\arg\min_{h\in\mathcal H_N}
\{\|\widehat b_{N,h} -\widehat b_{N,h_0}\|_{\delta}^{2} +\textrm{pen}(h)\}
\end{equation}
with
\begin{displaymath}
\textrm{pen}(h) :=
\frac{2}{(T - t_0)^2N^2}
\sum_{i = 1}^{N}\left\langle\int_{t_0}^{T}K_h(F(X_{s}^{i}) - F(.))dX_{s}^{i},
\int_{t_0}^{T}K_{h_0}(F(X_{s}^{i}) - F(.))dX_{s}^{i}\right\rangle_{\delta}
\textrm{ $;$ }
\forall h\in\mathcal H_N.
\end{displaymath}
The following theorem deals with a risk bound on the adaptive estimator
\begin{displaymath}
\widehat b_N(x) :=
\widehat b_{N,\widehat h}(x)
\textrm{ $;$ }x\in\mathbb R.
\end{displaymath}
%


%
\begin{theorem}\label{risk_bound_PCO_estimator_b}
Under Assumptions \ref{assumption_b_sigma}, \ref{assumption_K} and \ref{assumption_K_sigma_PCO}, if $t_0 > 0$, then there exist two deterministic constants $\mathfrak c_{\ref{risk_bound_PCO_estimator_b},1},\mathfrak c_{\ref{risk_bound_PCO_estimator_b},2} > 0$, not depending on $N$, such that for every $\vartheta\in (0,1)$ and $\lambda > 0$, with probability larger than $1 -\mathfrak c_{\ref{risk_bound_PCO_estimator_b},1}|\mathcal H_N|e^{-\lambda}$,
\begin{displaymath}
\|\widehat b_N - b\|_{\delta}^{2}
\leqslant
(1 +\vartheta)\min_{h\in\mathcal H_N}
\|\widehat b_{N,h} - b\|_{\delta}^{2} +
\frac{\mathfrak c_{\ref{risk_bound_PCO_estimator_b},2}}{\vartheta}
\left[\|b_{h_0} - b\|_{\delta}^{2} +\frac{(1 +\lambda)^3}{N}\right].
\end{displaymath}
\end{theorem}
\noindent
Theorem \ref{risk_bound_PCO_estimator_b} allows to establish a risk bound on our adaptive warped estimator for the usual norm $\|.\|_{2,A,B}$ on $\mathbb L^2([A,B],dx)$, where $A,B\in\mathbb R$ satisfy $-\Delta < A < B <\Delta$. First, since $\delta(I)\subset (0,\infty)$ for every compact subset $I$ of ${\rm supp}(\delta)^{\circ} = (-\Delta,\Delta)$, there exists $\delta_0 > 0$ such that $\delta(x) >\delta_0$ for every $x\in [A,B]$. Then, for any $\vartheta\in (0,1)$ and $\lambda > 0$, with probability larger than $1 -\mathfrak c_{\ref{risk_bound_PCO_estimator_b},1}|\mathcal H_N|e^{-\lambda}$,
\begin{eqnarray*}
 \|\widehat b_N - b\|_{2,A,B}^{2}
 & \leqslant &
 \frac{1}{\delta_0}\|\widehat b_N - b\|_{\delta}^{2}\\
 & \leqslant &
 \frac{1}{\delta_0}\left[
 (1 +\vartheta)\min_{h\in\mathcal H_N}
 \|\widehat b_{N,h} - b\|_{\delta}^{2} +
 \frac{\mathfrak c_{\ref{risk_bound_PCO_estimator_b},2}}{\vartheta}
 \left[\|b_{h_0} - b\|_{\delta}^{2} +\frac{(1 +\lambda)^3}{N}\right]\right].
\end{eqnarray*}
Moreover, since ${\rm supp}(\delta) = [-\Delta,\Delta]$ and $f(\mathbb R)\subset (0,\infty)$, there exists $f_1 > 0$ such that $f(x)\geqslant f_1$ for every $x\in {\rm supp}(\delta)$, and then
\begin{eqnarray*}
 & &
 (1 +\vartheta)\min_{h\in\mathcal H_N}
 \|\widehat b_{N,h} - b\|_{\delta}^{2} +
 \frac{\mathfrak c_{\ref{risk_bound_PCO_estimator_b},2}}{\vartheta}
 \left[\|b_{h_0} - b\|_{\delta}^{2} +\frac{(1 +\lambda)^3}{N}\right]\\
 & &
 \hspace{2cm}\leqslant
 \frac{1\vee\|\delta\|_{\infty}}{1\wedge f_1}\left[(1 +\vartheta)\min_{h\in\mathcal H_N}
 \|\widehat b_{N,h} - b\|_{f,-\Delta,\Delta}^{2} +
 \frac{\mathfrak c_{\ref{risk_bound_PCO_estimator_b},2}}{\vartheta}
 \left[\|b_{h_0} - b\|_{f,-\Delta,\Delta}^{2} +\frac{(1 +\lambda)^3}{N}\right]\right].
\end{eqnarray*}
Therefore, with probability larger than $1 -\mathfrak c_{\ref{risk_bound_PCO_estimator_b},1}|\mathcal H_N|e^{-\lambda}$,
\begin{displaymath}
\|\widehat b_N - b\|_{2,A,B}^{2}
\leqslant
\frac{1\vee\|\delta\|_{\infty}}{\delta_0(1\wedge f_1)}
\left[(1 +\vartheta)\min_{h\in\mathcal H_N}
\|\widehat b_{N,h} - b\|_{f,-\Delta,\Delta}^{2} +
\frac{\mathfrak c_{\ref{risk_bound_PCO_estimator_b},2}}{\vartheta}
\left[\|b_{h_0} - b\|_{f,-\Delta,\Delta}^{2} +\frac{(1 +\lambda)^3}{N}\right]\right]
\end{displaymath}
and, by Proposition \ref{risk_bound_warped_estimator}, there exist two constants $\mathfrak c_1,\mathfrak c_2 > 0$, not depending on $N$, such that
\begin{displaymath}
\mathbb E\left(\min_{h\in\mathcal H_N}
\|\widehat b_{N,h} - b\|_{f,-\Delta,\Delta}^{2}\right)
\leqslant
\mathfrak c_1\min_{h\in\mathcal H_N}\left\{h^2 +\frac{1}{Nh}\right\}
\quad {\rm and}\quad
\|b_{h_0} - b\|_{f,-\Delta,\Delta}^{2}\leqslant\mathfrak c_2h_{0}^{2}.
\end{displaymath}
%


%
\section{Numerical experiments}\label{section_numerical_experiments}
\begin{figure}[!ht]
\centering
\includegraphics[scale=0.45]{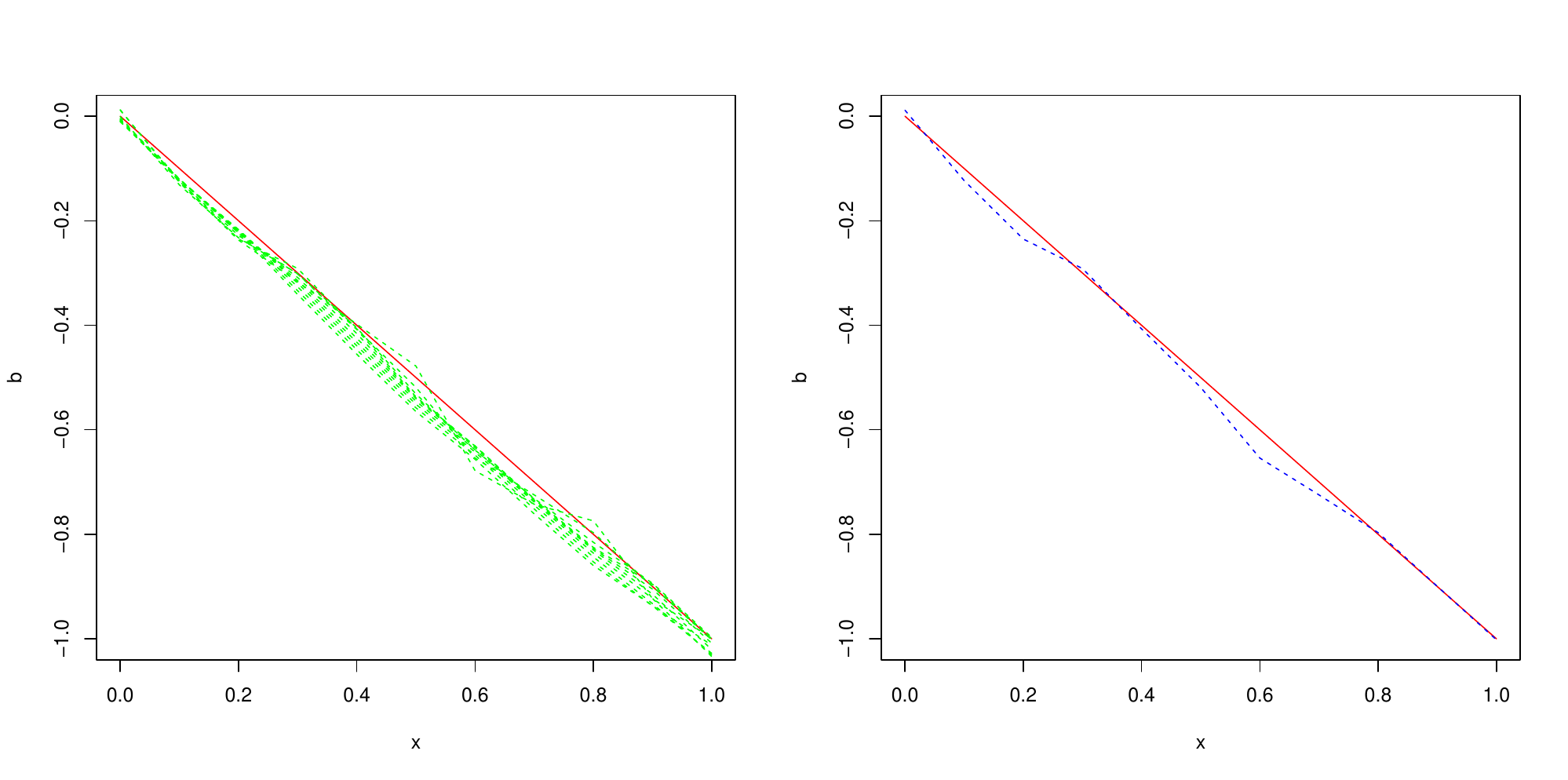} 
\caption{PCO adaptative estimator for Model 1 (Langevin equation), $\widehat h = 0.04$ and $h_{\rm oracle} = 0.02$.}
\label{plot_model_1}
\end{figure}

\begin{figure}[!h]
\centering
\includegraphics[scale=0.45]{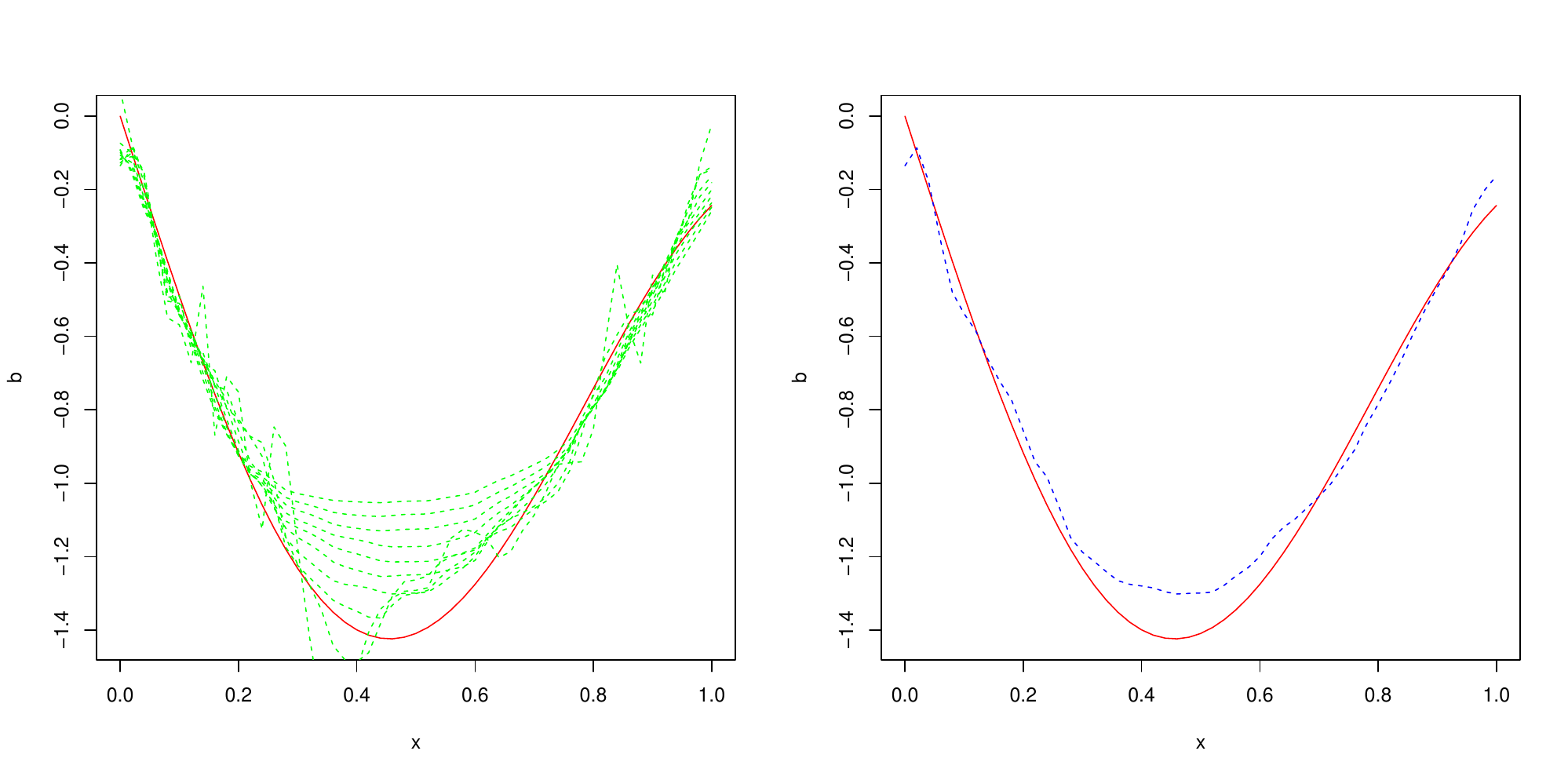} 
\caption{PCO adaptative estimator for Model 2, $\widehat h = 0.05$ and $h_{\rm oracle} = 0.07$.}
\label{plot_model_2}
\end{figure}

First of all, recall that since the function $F$ is unknown in practice, $\widehat b_{N,h}(x)$ is approximated by
\begin{displaymath}
\widetilde b_{N,h}(x) =
\widehat\beta_{N,h}(\widehat F_N;\widehat F_N(x))
\quad {\rm with}\quad
\widehat F_N(x) =
\frac{1}{N(T - t_0)}\sum_{i = 1}^{N}\int_{t_0}^{T}\mathbf 1_{X_{t}^{i}\leqslant x}dt.
\end{displaymath}
A discrete-time approximate version of this estimator is computed on datasets generated by two different SDEs. In each case, the bandwidth of our estimator is selected via the PCO method introduced in Subsection \ref{section_PCO}. The first model is the so-called Langevin equation, and the second one is a nonlinear SDE with a multiplicative noise:
\begin{itemize}
 \item[1.] $\displaystyle{X_t = x_0 -\int_{0}^{t}X_sds + 0.1\cdot W_t}$ ; $t\in [0,T]$.
 \item[2.] $\displaystyle{X_t = x_0 -\int_{0}^{t}(X_s + \sin(4X_s))ds + 0.1\int_{0}^{t}(2 +\cos(X_s))dW_s}$ ; $t\in [0,T]$.
\end{itemize}
The models and the estimator are implemented by taking $N = 100$, $n = 50$, $T = 5$, $x_0 = 2$, $t_0 = 0$ and $K =\delta =\rho$, where $\rho$ is the kernel introduced in Subsection \ref{section_PCO}. For Model 1, the estimator of the drift function is computed for the bandwidths set
\begin{displaymath}
\mathcal H_1 :=\{0.02k\textrm{ $;$ }k = 1,\dots,10\},
\end{displaymath}
and for Model 2, it is computed for the bandwidths set
\begin{displaymath}
\mathcal H_2 :=\{0.01k\textrm{ $;$ }k = 1,\dots,10\}.
\end{displaymath}
For each of the previous models, on Figures \ref{plot_model_1} and \ref{plot_model_2} respectively, the true drift function (in red) and the PCO adaptive estimator (in blue) are plotted on the left-hand side, and the beam of proposals is plotted in green on the right-hand side. On Figure \ref{plot_model_1}, one can see that the drift function is well estimated by the PCO adaptive estimator, with a MSE equal to $7.12\cdot 10^{-4}$ compared to a MSE equal to $5.11\cdot 10^{-4}$ for the oracle estimator. On Figure \ref{plot_model_2}, one can see that the drift function of Model 2 is still well estimated by our PCO adaptive estimator. However, note that there is a degradation of the MSE, which is equal to $7.47\cdot 10^{-3}$ for the adaptative estimator and equal to $4.99\cdot 10^{-3}$ for the oracle estimator. This is probably related to both the nonlinearity of the drift function and the multiplicative noise. In some cases, especially for Model 2, when the bandwidth is too small, the estimation degrades, but the PCO method selects a higher value of $h$ which performs better.
Finally, for each model, Table \ref{table_MSE} gathers the mean MSE of 100 PCO estimations of the drift function as well as the mean MSE of the corresponding 100 oracle estimations. The mean MSEs are globally low, but higher for Model 2 with a nonlinear drift function than for Model 1 with a linear one. However, regardless of the complexity of the model, one can notice that for both models, the mean MSE of the PCO estimations remains relatively close to the mean MSE of the corresponding oracle estimations. This means that our PCO method performs well in practice.
\begin{table}[h!]
\begin{center}
\begin{tabular}{|c|c|c|}
\hline
 & PCO & Oracle \\
 \hline\hline
 Model 1 & $8.30\cdot 10^{-4}$ & $6.28\cdot 10^{-4}$\\
 \hline
 Model 2 & $5.86\cdot 10^{-3}$ & $4.77\cdot 10^{-3}$\\
 \hline
\end{tabular}
\medskip
\caption{Mean MSEs of 100 PCO adaptive estimations compared to the oracle estimations.}\label{table_MSE}
\end{center}
\end{table}
%


%
\section{Proofs}\label{section_proofs}
%


%
\subsection{Proof of Proposition \ref{risk_bound_warped_estimator}}
For the sake of readability, without loss of generality, Proposition \ref{risk_bound_warped_estimator} is proved for $t_0 = 0$. First of all,
\begin{eqnarray*}
\mathbb E(\|\widehat b_{N,h} - b\|_{f,A,B}^{2})\leqslant
\int_{A}^{B}\mathfrak b(x)^2f(x)dx
+\int_{-\infty}^{\infty}\mathfrak v(x)f(x)dx
\end{eqnarray*}
where $\mathfrak b(x)$ (resp. $\mathfrak v(x)$) is the bias (resp. the variance) term of $\widehat b_{N,h}(x)$ for every $x\in [A,B]$ (resp. $x\in\mathbb R$). On the one hand, let us find a suitable bound on the $f$-weighted integrated variance of our warped kernel estimator. For any $x\in\mathbb R$, since $X^1,\dots,X^N$ are independent copies of $X$,
\begin{eqnarray*}
 \mathfrak v(x)
 & = &
 \textrm{var}\left(\frac{1}{NT}
 \sum_{i = 1}^{N}\int_{0}^{T}K_h(F(X_{t}^{i}) - F(x))dX_{t}^{i}\right)
 \leqslant
 \frac{1}{NT^2}
 \mathbb E\left[\left(\int_{0}^{T}K_h(F(X_t) - F(x))dX_t\right)^2\right]\\
 & \leqslant &
 \frac{2}{N}\mathbb E\left[
 \left(\int_{0}^{T}K_h(F(X_t) - F(x))b(X_t)\frac{dt}{T}\right)^2 +
 \frac{1}{T^2}\left(\int_{0}^{T}K_h(F(X_t) - F(x))\sigma(X_t)dW_t\right)^2\right].
\end{eqnarray*}
In the right-hand side of the previous inequality, Jensen's inequality on the first term and the isometry property for It\^o's integral on the second one give
\begin{eqnarray*}
 \mathfrak v(x) & \leqslant &
 \frac{2}{NT}\int_{0}^{T}\mathbb E[K_h(F(X_t) - F(x))^2b(X_t)^2]dt +
 \frac{2}{NT^2}\int_{0}^{T}\mathbb E[K_h(F(X_t) - F(x))^2\sigma(X_t)^2]dt\\
 & = &
 \frac{2}{N}\int_{-\infty}^{\infty}K_h(F(z) - F(x))^2b(z)^2f(z)dz +
 \frac{2}{NT}\int_{-\infty}^{\infty}K_h(F(z) - F(x))^2\sigma(z)^2f(z)dz.
\end{eqnarray*}
Since $K$ is symmetric, $K\in\mathbb L^2(\mathbb R,dx)$, $b,\sigma\in\mathbb L^2(\mathbb R,f(x)dx)$ and $F$ is one-to-one from $\mathbb R$ into $(0,1)$ (see Assumption \ref{assumption_K} and Remark \ref{F_properties}),
\begin{eqnarray*}
 \int_{-\infty}^{\infty}\mathfrak v(x)f(x)dx
 & \leqslant &
 \frac{2}{N}\int_{\mathbb R^2}
 K_h(F(z) - F(x))^2
 \left[b(z)^2 +\frac{1}{T}\sigma(z)^2\right]f(z)f(x)dxdz\\
 & = &
 \frac{2}{N}\int_{-\infty}^{\infty}
 \left[b(z)^2 +\frac{1}{T}\sigma(z)^2\right]f(z)\\
 & &
 \hspace{2cm}\times
 \int_{F(\mathbb R)}K_h(F(z) - y)^2f(F^{-1}(y))\frac{dy}{F'(F^{-1}(y))}dz\\
 & = &
 \frac{2}{Nh}\int_{-\infty}^{\infty}
 \left[b(z)^2 +\frac{1}{T}\sigma(z)^2\right]f(z)
 \int_{-F(z)/h}^{(1 - F(z))/h}K(y)^2dydz\\
 & \leqslant &
 \frac{2\|K\|_{2}^{2}}{Nh}\left(\|b\|_{f}^{2} +
 \frac{1}{T}\|\sigma\|_{f}^{2}\right).
\end{eqnarray*}
On the other hand, let us find a suitable bound on the $[A,B]$-$f$-weighted integrated squared-bias of $\widehat b_{N,h}$. For every $x\in [A,B]$, since $X^1,\dots,X^N$ are independent copies of $X$, and since It\^o's integral restricted to $\mathbb H^2$ is a martingale-valued map,
\begin{eqnarray*}
 \mathfrak b(x) + b(x)
 & = &
 \mathbb E\left[\frac{1}{NT}
 \sum_{i = 1}^{N}\int_{0}^{T}K_h(F(X_{t}^{i}) - F(x))dX_{t}^{i}\right] =
 \frac{1}{T}
 \mathbb E\left(\int_{0}^{T}K_h(F(X_t) - F(x))dX_t\right)\\
 & = &
 \frac{1}{T}\left[\mathbb E\left(\int_{0}^{T}K_h(F(X_t) - F(x))b(X_t)dt\right)
 +\mathbb E\left(\int_{0}^{T}K_h(F(X_t) - F(x))\sigma(X_t)dW_t\right)\right]\\
 & = &
 \frac{1}{T}\int_{0}^{T}\mathbb E(K_h(F(X_t) - F(x))b(X_t))dt\\
 & &
 \hspace{1cm}
 =\int_{-\infty}^{\infty}K_h(F(z) - F(x))b(z)f(z)dz\\
 & &
 \hspace{2cm}
 =\int_{0}^{1}K_h(y - F(x))b(F^{-1}(y))dy
 = [K_h\ast ((b\circ F^{-1})\mathbf 1_{(0,1)})](F(x)).
\end{eqnarray*}
Moreover, since $b\in\mathbb L^2(\mathbb R,f(x)dx)$ by Remark \ref{F_properties},
\begin{eqnarray*}
 \|[K_h\ast ((b\circ F^{-1})\mathbf 1_{(0,1)})](F(.))\|_{f,A,B}^{2} & \leqslant &
 \int_{-\infty}^{\infty}[K_h\ast ((b\circ F^{-1})\mathbf 1_{(0,1)})](F(x))^2f(x)dx\\
 & = &
 \int_{0}^{1}[K_h\ast ((b\circ F^{-1})\mathbf 1_{(0,1)})](y)^2dy\\
 & \leqslant &
 \|K\|_{1}^{2}\int_{0}^{1}b(F^{-1}(y))^2dy =
 \|K\|_{1}^{2}\|b\|_{f}^{2} <\infty.
\end{eqnarray*}
Therefore,
\begin{displaymath}
\int_{A}^{B}\mathfrak b(x)^2f(x)dx =
\|[K_h\ast ((b\circ F^{-1})\mathbf 1_{(0,1)})](F(.)) - b\|_{f,A,B}^{2}.
\end{displaymath}
This concludes the proof.
%


%
\subsection{Proof of Theorem \ref{risk_bound_PCO_estimator_b}}\label{section_proof_PCO}
The proof of Theorem \ref{risk_bound_PCO_estimator_b} relies on the following technical results.
%


%
\begin{proposition}\label{kernel_properties_Phi}
Under Assumptions \ref{assumption_b_sigma}, \ref{assumption_K} and \ref{assumption_K_sigma_PCO},
\begin{displaymath}
\frac{1}{T - t_0}\int_{t_0}^{T}K_h(F(X_t) - F(x))dX_t =\Phi_h(X,x)
\textrm{ $;$ }
\forall x\in\mathbb R\textrm{, }
\forall h > 0,
\end{displaymath}
where $(x,h,\varphi)\mapsto\Phi_h(\varphi,x)$ is the map from $\mathbb R\times (0,\infty)\times C^0([0,T];\mathbb R)$ into $\mathbb R$ defined by
\begin{displaymath}
\Phi_h(\varphi,x) :=
\frac{1}{T - t_0}\left[\int_{\varphi(t_0)}^{\varphi(T)}K_h(F(z) - F(x))dz -
\frac{1}{2h^2}\int_{t_0}^{T}K'\left(\frac{F(\varphi(t)) - F(x)}{h}\right)\sigma(\varphi(t))^2f(\varphi(t))dt
\right]
\end{displaymath}
for every $x\in\mathbb R$, $h > 0$ and $\varphi\in C^0([t_0,T];\mathbb R)$. Moreover,
\begin{enumerate}
 \item There exists a constant $\mathfrak c_{\ref{kernel_properties_Phi},1} > 0$ such that, for every $h\in (0,\Delta_0]$ and $\varphi\in C^0([t_0,T];\mathbb R)$,
 \begin{displaymath}
 \|\Phi_h(\varphi,.)\|_{\delta}^{2}
 \leqslant
 \frac{\mathfrak c_{\ref{kernel_properties_Phi},1}}{h^3}.
 \end{displaymath}
 \item There exists a constant $\mathfrak c_{\ref{kernel_properties_Phi},2} > 0$ such that, for every $h,h'\in (0,\Delta_0]$,
 \begin{displaymath}
 \mathbb E(\langle\Phi_h(X^1,.),\Phi_{h'}(X^2,.)\rangle_{\delta}^{2})
 \leqslant
 \mathfrak c_{\ref{kernel_properties_Phi},2}\mathfrak m(h')
 \end{displaymath}
 with
 \begin{displaymath}
 \mathfrak m(h') =\mathbb E(\|\Phi_{h'}(X,.)\|_{\delta}^{2}).
 \end{displaymath}
 \item There exists a constant $\mathfrak c_{\ref{kernel_properties_Phi},3} > 0$ such that, for every $h\in (0,\Delta_0]$ and $\varphi\in\mathbb L^2(\mathbb R,dx)$,
 \begin{displaymath}
 \mathbb E(\langle\Phi_h(X,.),\varphi\rangle_{\delta}^{2})
 \leqslant\mathfrak c_{\ref{kernel_properties_Phi},3}
 \|\varphi\|_{\delta}^{2}.
 \end{displaymath}
 \item If $t_0 > 0$, then there exists a deterministic constant $\mathfrak c_{\ref{kernel_properties_Phi},4} > 0$ such that, for every $h,h'\in\mathcal H_N$,
 \begin{displaymath}
 |\langle\Phi_h(X,.),b_{h'}\rangle_{\delta}|
 \leqslant\mathfrak c_{\ref{kernel_properties_Phi},4}
 \quad
 {\rm a.s.}
 \end{displaymath}
\end{enumerate}
\end{proposition}
%


%
\begin{lemma}\label{bound_U_statistics}
Consider
\begin{equation}\label{U_statistic_definition}
U_{h,h'}(N) :=\sum_{i\not= j}
\langle\Phi_h(X^i,.) - b_h,
\Phi_{h'}(X^j,.) - b_{h'}\rangle_{\delta}
\textrm{ $;$ }
\forall h,h'\in\mathcal H_N.
\end{equation}
Under Assumptions \ref{assumption_b_sigma}, \ref{assumption_K} and \ref{assumption_K_sigma_PCO}, there exists a deterministic constant $\mathfrak c_{\ref{bound_U_statistics}} > 0$, not depending on $N$, such that for every $\theta\in (0,1)$ and $\lambda > 0$, with probability larger than $1 - 5.4|\mathcal H_N|e^{-\lambda}$,
\begin{eqnarray*}
 & & \sup_{h\in\mathcal H_N}
 \left\{\frac{|U_{h,h_0}(N)|}{N^2}
 -\frac{\theta\mathfrak m(h)}{N}\right\}
 \leqslant
 \frac{\mathfrak c_{\ref{bound_U_statistics}}(1 +\lambda)^3}{\theta N}\\
 & &
 \hspace{4cm}\textrm{and }
 \sup_{h\in\mathcal H_N}
 \left\{\frac{|U_{h,h}(N)|}{N^2}
 -\frac{\theta\mathfrak m(h)}{N}\right\}
 \leqslant
 \frac{\mathfrak c_{\ref{bound_U_statistics}}(1 +\lambda)^3}{\theta N}.
\end{eqnarray*}
\end{lemma}
%


%
\begin{lemma}\label{bound_trace_term}
Consider
\begin{displaymath}
V_h(N) :=\frac{1}{N}\sum_{i = 1}^{N}\|\Phi_h(X^i,.) - b_h\|_{\delta}^{2}
\textrm{ $;$ }
\forall h\in\mathcal H_N.
\end{displaymath}
Under Assumptions \ref{assumption_b_sigma}, \ref{assumption_K} and \ref{assumption_K_sigma_PCO}, there exists a deterministic constant $\mathfrak c_{\ref{bound_trace_term}} > 0$, not depending on $N$, such that for every $\theta\in (0,1)$ and $\lambda > 0$, with probability larger than $1 - 2|\mathcal H_N|e^{-\lambda}$,
\begin{displaymath}
\sup_{h\in\mathcal H_N}\left\{
\frac{1}{N}|V_h(N) -\mathfrak m(h)| -\frac{\theta\mathfrak m(h)}{N}\right\}
\leqslant
\frac{\mathfrak c_{\ref{bound_trace_term}}(1 +\lambda)}{\theta N}.
\end{displaymath}
\end{lemma}
%


%
\begin{lemma}\label{bound_crossed_term}
Consider
\begin{equation}\label{W_statistic_definition}
W_{h,h'}(N) :=
\langle\widehat b_{N,h} - b_h,b_{h'} - b\rangle_{\delta}
\textrm{ $;$ }
\forall h,h'\in\mathcal H_N.
\end{equation}
Under Assumptions \ref{assumption_b_sigma}, \ref{assumption_K} and \ref{assumption_K_sigma_PCO}, there exists a deterministic constant $\mathfrak c_{\ref{bound_crossed_term}} > 0$, not depending on $N$, such that for every $\theta\in (0,1)$ and $\lambda > 0$, with probability larger than $1 - 2|\mathcal H_N|e^{-\lambda}$,
\begin{eqnarray*}
 & & \sup_{h\in\mathcal H_N}\{
 |W_{h,h_0}(N)| -\theta\|b_{h_0} - b\|_{\delta}^{2}\}
 \leqslant
 \frac{\mathfrak c_{\ref{bound_crossed_term}}(1 +\lambda)^2}{\theta N},\\
 & &
 \hspace{2cm}\sup_{h\in\mathcal H_N}\{
 |W_{h_0,h}(N)| -\theta\|b_h - b\|_{\delta}^{2}\}
 \leqslant
 \frac{\mathfrak c_{\ref{bound_crossed_term}}(1 +\lambda)^2}{\theta N}\\
 & &
 \hspace{4cm}\textrm{and}
 \quad\sup_{h\in\mathcal H_N}\{
 |W_{h,h}(N)| -\theta\|b_h - b\|_{\delta}^{2}\}
 \leqslant
 \frac{\mathfrak c_{\ref{bound_crossed_term}}(1 +\lambda)^2}{\theta N}.
\end{eqnarray*}
\end{lemma}
\noindent
The proof of Proposition \ref{kernel_properties_Phi}, which is specific to our warped kernel estimator of $b$, is postponed to Subsubsection \ref{Section_proof_kernel_properties_Phi}. By using our Proposition \ref{kernel_properties_Phi} instead of Marie and Rosier \cite{MR23}, Lemma 3, the proofs of Lemmas \ref{bound_U_statistics}, \ref{bound_trace_term} and \ref{bound_crossed_term} remain the same as those of Marie and Rosier \cite{MR23}, Lemmas 4, 5 and 6 respectively.
%


%
\subsubsection{Steps of the proof}
The proof of Theorem \ref{risk_bound_PCO_estimator_b}.(1) is dissected in four steps.
\\
\\
\textbf{Step 1.} This first step provides a suitable decomposition of $\|\widehat b_{N,\widehat h} - b\|_{\delta}^{2}$. First,
\begin{displaymath}
\|\widehat b_{N,\widehat h} - b\|_{\delta}^{2} =
\|\widehat b_{N,\widehat h} -\widehat b_{N,h_0}\|_{\delta}^{2} +
\|\widehat b_{N,h_0} - b\|_{\delta}^{2} -
2\langle\widehat b_{N,h_0} -\widehat b_{N,\widehat h},
\widehat b_{N,h_0} - b\rangle_{\delta}.
\end{displaymath}
Then, by (\ref{penalty_proposal}), for any $h\in\mathcal H_N$, 
\begin{eqnarray}
 \|\widehat b_{N,\widehat h} - b\|_{\delta}^{2}
 & \leqslant &
 \|\widehat b_{N,h} -\widehat b_{N,h_0}\|_{\delta}^{2} +
 \textrm{pen}(h) -\textrm{pen}(\widehat h)
 \nonumber\\
 & &
 \hspace{2cm}
 +\|\widehat b_{N,h_0} - b\|_{\delta}^{2} -
 2\langle\widehat b_{N,h_0} -\widehat b_{N,\widehat h},
 \widehat b_{N,h_0} - b\rangle_{\delta}
 \nonumber\\
 & \leqslant &
 \|\widehat b_{N,h} - b\|_{\delta}^{2} +
 \textrm{pen}(h) -\textrm{pen}(\widehat h) - 2\langle\widehat b_{N,h} -\widehat b_{N,\widehat h},
 \widehat b_{N,h_0} - b\rangle_{\delta}
 \nonumber\\
 \label{risk_bound_PCO_estimator_b_1}
 & = &
 \|\widehat b_{N,h} - b\|_{\delta}^{2} -\psi_N(h) +\psi_N(\widehat h)
\end{eqnarray}
where 
\begin{displaymath}
\psi_N(h) :=
2\langle\widehat b_{N,h} - b,\widehat b_{N,h_0} -
b\rangle_{\delta} -\textrm{pen}(h).
\end{displaymath}
Let's complete the decomposition of $\|\widehat b_{N,\widehat h} - b\|_{\delta}^{2}$ by writing
\begin{displaymath}
\psi_N(h) = 2(\psi_{1,N}(h) + \psi_{2,N}(h) + \psi_{3,N}(h)),
\end{displaymath}
where
\begin{eqnarray*}
 \psi_{1,N}(h) & := &
 \frac{1}{(T - t_0)^2N^2}
 \sum_{i = 1}^{N}\left\langle\int_{t_0}^{T}K_h(F(X_{s}^{i}) - F(\cdot))dX_{s}^{i},
 \int_{t_0}^{T}K_{h_0}(F(X_{s}^{i}) - F(\cdot))dX_{s}^{i}\right\rangle_{\delta}\\
 & &
 \hspace{8cm}
 +\frac{U_{h,h_0}(N)}{N^2} -\frac{1}{2}{\rm pen}(h)
 =\frac{U_{h,h_0}(N)}{N^2},\\
 \psi_{2,N}(h) & := & -\frac{1}{N^2}\left(
 \sum_{i = 1}^{N}
 \left\langle\frac{1}{T - t_0}
 \int_{t_0}^{T}K_{h_0}(F(X_{s}^{i}) - F(\cdot))dX_{s}^{i},
 b_h\right\rangle_{\delta} +\right.\\
 & &
 \hspace{2cm}
 \left. +\sum_{i = 1}^{N}
 \left\langle\frac{1}{T - t_0}
 \int_{t_0}^{T}K_h(F(X_{s}^{i}) - F(\cdot))dX_{s}^{i},b_{h_0}\right\rangle_{\delta}\right)
 +\frac{1}{N}\langle b_{h_0},b_h\rangle_{\delta}\textrm{ and}\\
 \psi_{3,N}(h)
 & := &
 W_{h,h_0}(N) + W_{h_0,h}(N) +\langle b_h - b,b_{h_0} - b\rangle_{\delta}.
\end{eqnarray*}	
\textbf{Step 2.} This step deals with bounds on $\mathbb E(\psi_{j,N}(h))$ and $\mathbb E(\psi_{j,N}(\widehat h))$ for $j = 1,2,3$.
\begin{itemize}
 \item By Lemma \ref{bound_U_statistics}, for any $\lambda > 0$ and $\theta\in (0,1)$, with probability larger than $1 - 5.4|\mathcal H_N|e^{-\lambda}$,
 \begin{displaymath}
 |\psi_{1,N}(h)|
 \leqslant \frac{\theta\mathfrak m(h)}{N} +
 \frac{\mathfrak c_{\ref{bound_U_statistics}}(1 +\lambda)^3}{\theta N}
 \quad {\rm and}\quad
 |\psi_{1,N}(\widehat h)|
 \leqslant\frac{\theta\mathfrak m(\widehat h)}{N} +
 \frac{\mathfrak c_{\ref{bound_U_statistics}}(1 +\lambda)^3}{\theta N}.
 \end{displaymath}
 \item On the one hand, for any $h,h'\in\mathcal H_N$, consider 
 \begin{displaymath}
 \Psi_{2,N}(h,h') :=
 \frac{1}{N}\sum_{i = 1}^{N}
 \langle\Phi_h(X^i,.),b_{h'}\rangle_{\delta}.
 \end{displaymath}
 By Proposition \ref{kernel_properties_Phi}.(4),
 \begin{displaymath}
 |\Psi_{2,N}(h,h')|
 \leqslant
 \frac{1}{N}\sum_{i = 1}^{N}
 \left|\int_{-\infty}^{\infty}
 \Phi_h(X^i,x)b_{h'}(x)\delta(x)dx\right|\\
 \leqslant
 \mathfrak c_{\ref{kernel_properties_Phi},4}
 \quad {\rm a.s.}
 \end{displaymath}
 On the other hand, since
 \begin{eqnarray*}
  \|b_h\mathbf 1_{[-\Delta,\Delta]}\|_{\infty}
  & \leqslant &
  \|K\|_1\|b\|_{\infty,I_{\Delta,F}}\\
  & &
  {\rm and}\quad
  \|b_{h_0}\delta\|_1
  \leqslant
  \|b_{h_0}\mathbf 1_{[-\Delta,\Delta]}\|_{\infty}
  \|\delta\|_1\leqslant
  \|K\|_1\|b\|_{\infty,I_{\Delta,F}}
 \end{eqnarray*}
 by Inequality (\ref{kernel_properties_Phi_5}),
 \begin{displaymath}
 |\langle b_h,b_{h_0}\rangle_{\delta}| 
 \leqslant
 \|b_h\mathbf 1_{[-\Delta,\Delta]}\|_{\infty}
 \|b_{h_0}\delta\|_1
 \leqslant
 \|K\|_{1}^{2}\|b\|_{\infty,I_{\Delta,F}}^{2}.
 \end{displaymath}	
 Then, there exists a deterministic constant $\mathfrak c_1 > 0$, not depending on $N$ and $h$, such that
 \begin{displaymath}
 |\psi_{2,N}(h)|\leqslant\frac{\mathfrak c_1}{N}
 \quad {\rm and}\quad
 |\psi_{2,N}(\widehat h)|\leqslant
 \sup_{h'\in\mathcal H_N}
 |\psi_{2,N}(h')|\leqslant
 \frac{\mathfrak c_1}{N}
 \quad {\rm a.s.}
 \end{displaymath}
\item By Lemma \ref{bound_crossed_term} and Cauchy-Schwarz's inequality, with probability larger that $1 - |\mathcal H_N|e^{-\lambda}$,
 \begin{eqnarray*}
  |\psi_{3,N}(h)|
  & \leqslant &
  \frac{\theta}{4}(\|b_h - b\|_{\delta}^{2} +
  \|b_{h_0} - b\|_{\delta}^{2}) +
  \frac{8\mathfrak c_{\ref{bound_crossed_term}}(1 +\lambda)^2}{\theta N}\\
  & &
  \hspace{2cm}
  + 2\times\frac{1}{2^{1/2}}\left(\frac{\theta}{2}\right)^{1/2}
  \|b_h - b\|_{\delta}\times
  \frac{1}{2^{1/2}}\left(\frac{2}{\theta}\right)^{1/2}\|b_{h_0} - b\|_{\delta}\\
  & \leqslant &
  \frac{\theta}{2}\|b_h - b\|_{\delta}^{2} +
  \left(\frac{\theta}{4} +\frac{1}{\theta}\right)
  \|b_{h_0} - b\|_{\delta}^{2} +
  \frac{8\mathfrak c_{\ref{bound_crossed_term}}(1 +\lambda)^2}{\theta N}
 \end{eqnarray*}	
 and
 \begin{displaymath}
 |\psi_{3,N}(\widehat h)|\leqslant
 \frac{\theta}{2}\|b_{\widehat h} - b\|_{\delta}^{2} +
 \left(\frac{\theta}{4} +\frac{1}{\theta}\right)
 \|b_{h_0} - b\|_{\delta}^{2} +
 \frac{8\mathfrak c_{\ref{bound_crossed_term}}(1 +\lambda)^2}{\theta N}.
 \end{displaymath}
\end{itemize}
\textbf{Step 3.} Let us establish that there exist two deterministic constants $\mathfrak c_2,\overline{\mathfrak c}_2 > 0$, not depending on $N$ and $\theta$, such that with probability larger than $1 -\overline{\mathfrak c}_2|\mathcal H_N|e^{-\lambda}$,
\begin{displaymath}
\sup_{h\in\mathcal H_N}\left\{
\|\widehat b_{N,h} - b\|_{\delta}^{2} -
(1 +\theta)\left(\|b_h - b\|_{\delta}^{2} +
\frac{\mathfrak m(h)}{N}\right)\right\}
\leqslant
\frac{\mathfrak c_2(1 +\lambda)^3}{\theta N}
\end{displaymath}
and
\begin{displaymath}
\sup_{h\in\mathcal H_N}\left\{
\|b_h - b\|_{\delta}^{2} +\frac{\mathfrak m(h)}{N} -
\frac{1}{1 -\theta}\|\widehat b_{N,h} - b\|_{\delta}^{2}\right\}
\leqslant
\frac{\mathfrak c_2(1 +\lambda)^3}{\theta(1 -\theta)N}.
\end{displaymath}
On the one hand, note that
\begin{displaymath}
\|\widehat b_{N,h} - b\|_{\delta}^{2} -
(1 +\theta)\left(\|b_h - b\|_{\delta}^{2} +\frac{\mathfrak m(h)}{N}\right)
\end{displaymath}
can be written
\begin{displaymath}
\|\widehat b_{N,h} - b_h\|_{\delta}^{2} -\frac{(1 +\theta)\mathfrak m(h)}{N} +
2W_h(N) -\theta\|b_h - b\|_{\delta}^{2},
\end{displaymath}
where $W_h(N) := W_{h,h}(N)$ (see (\ref{W_statistic_definition})). Moreover, for any $h\in\mathcal H_N$,
\begin{equation}\label{risk_bound_PCO_estimator_b_2}
\|\widehat b_{N,h} - b_h\|_{\delta}^{2} =
\frac{U_h(N)}{N^2} +\frac{V_h(N)}{N}
\end{equation}
with $U_h(N) = U_{h,h}(N)$ (see (\ref{U_statistic_definition})). So, with probability larger than $1 -\overline{\mathfrak c}_2|\mathcal H_N|e^{-\lambda}$,
\begin{displaymath}
\sup_{h\in\mathcal H_N}\left\{
\|\widehat b_{N,h} - b_h\|_{\delta}^{2} -
\frac{(1 +\theta)\mathfrak m(h)}{N}\right\}
\leqslant
\frac{2(\mathfrak c_{\ref{bound_U_statistics}} +\mathfrak c_{\ref{bound_trace_term}})(1 +\lambda)^3}{\theta N}
\end{displaymath}
by Lemmas \ref{bound_U_statistics} and \ref{bound_trace_term}, and then
\begin{displaymath}
\sup_{h\in\mathcal H_N}\left\{
\|\widehat b_{N,h} - b\|_{\delta}^{2} - (1 +\theta)\left(\|b_h - b\|_{\delta}^{2} +\frac{\mathfrak m(h)}{N}\right)\right\}
\leqslant\frac{\mathfrak c_2(1 +\lambda)^3}{\theta N}
\end{displaymath}
by Lemma \ref{bound_crossed_term}. On the other hand, for any $h\in\mathcal H_N$,
\begin{displaymath}
\|b_h - b\|_{\delta}^{2} =
\|\widehat b_{N,h} - b\|_{\delta}^{2}
-\|\widehat b_{N,h} - b_h\|_{\delta}^{2} - 2W_h(N).
\end{displaymath}
Then,
\begin{displaymath}
(1 -\theta)\left(\|b_h - b\|_{\delta}^{2} +\frac{\mathfrak m(h)}{N}\right)
-\|\widehat b_{N,h} - b\|_{\delta}^{2}
\leqslant
2|W_h(N)| -\theta\|b_h - b\|_{\delta}^{2} +
\Lambda_h(N) -\frac{\theta\mathfrak m(h)}{N}
\end{displaymath}
where
\begin{displaymath}
\Lambda_h(N) :=\left|
\|\widehat b_{N,h} - b_h\|_{\delta}^{2} -
\frac{\mathfrak m(h)}{N}\right|.
\end{displaymath}
By Equality (\ref{risk_bound_PCO_estimator_b_2}),
\begin{displaymath}
\Lambda_h(N) =
\left|\frac{U_h(N)}{N^2} +\frac{V_h(N)}{N} -\frac{\mathfrak m(h)}{N}\right|.
\end{displaymath}
By Lemmas \ref{bound_trace_term} and \ref{bound_U_statistics}, there exist two deterministic constants $\mathfrak c_3,\overline{\mathfrak c}_3 > 0$, not depending $N$ and $\theta$, such that with probability larger than $1 -\overline{\mathfrak c}_3|\mathcal H_N|e^{-\lambda}$,
\begin{displaymath}
\sup_{h\in\mathcal H_N}\left\{
\Lambda_h(N) -\theta\frac{\mathfrak m(h)}{N}\right\}
\leqslant
\frac{\mathfrak c_3(1 +\lambda)^3}{\theta N}.
\end{displaymath}
By Lemma \ref{bound_crossed_term}, with probability larger than $1 - 2|\mathcal H_N|e^{-\lambda}$,
\begin{displaymath}
\sup_{h\in\mathcal H_N}\{2|W_h(N)| -\theta\|b_h - b\|_{\delta}^{2}
\}\leqslant
\frac{4\mathfrak c_{\ref{bound_crossed_term}}(1 +\lambda)^2}{\theta N}.
\end{displaymath}
Therefore, with probability larger than $1 -\overline{\mathfrak c}_2|\mathcal H_N|e^{-\lambda}$,
\begin{displaymath}
\sup_{h\in\mathcal H_N}\left\{
\|b_h - b\|_{\delta}^{2} +\frac{\mathfrak m(h)}{N} -
\frac{1}{1 -\theta}\|\widehat b_{N,h} - b\|_{\delta}^{2}\right\}
\leqslant
\frac{\mathfrak c_2(1 +\lambda)^3}{\theta(1 -\theta)N}.
\end{displaymath}
\textbf{Step 4.} By Step 2, there exist two deterministic constants $\mathfrak c_4,\overline{\mathfrak c}_4 > 0$, not depending on $N$, $\theta$, $h$ and $h_0$, such that with probability larger than $1 -\overline{\mathfrak c}_4|\mathcal H_N|e^{-\lambda}$,
\begin{displaymath}
|\psi_N(h)|
\leqslant\theta\left(\|b_h - b\|_{\delta}^{2} +\frac{\mathfrak m(h)}{N}\right) +
\left(\frac{\theta}{2} +\frac{2}{\theta}\right)\|b_{h_0} - b\|_{\delta}^{2} +
\frac{\mathfrak c_4(1 +\lambda)^3}{\theta N}
\end{displaymath} 
and
\begin{displaymath}
|\psi_N(\widehat h)|
\leqslant\theta\left(\|b_{\widehat h} - b\|_{\delta}^{2} +
\frac{\mathfrak m(\widehat h)}{N}\right) +
\left(\frac{\theta}{2} +\frac{2}{\theta}\right)\|b_{h_0} - b\|_{\delta}^{2} +
\frac{\mathfrak c_4(1 +\lambda)^3}{\theta N}.
\end{displaymath} 
Then, by Step 3, there exist two deterministic constants $\mathfrak c_5,\overline{\mathfrak c}_5 > 0$, not depending on $N$, $\theta$, $h$ and $h_0$, such that with probability larger than $1 -\overline{\mathfrak c}_5|\mathcal H_N|e^{-\lambda}$,
\begin{displaymath}
|\psi_N(h)|
\leqslant
\frac{\theta}{1 -\theta}
\|\widehat b_{N,h} - b\|_{\delta}^{2} +
\left(\frac{\theta}{2} +\frac{2}{\theta}\right)\|b_{h_0} - b\|_{\delta}^{2} +
\mathfrak c_5\left(\frac{1}{\theta} +\frac{1}{1 -\theta}\right)
\frac{(1 +\lambda)^3}{N}
\end{displaymath}	
and
\begin{displaymath}
|\psi_N(\widehat h)|
\leqslant
\frac{\theta}{1 -\theta}
\|\widehat b_{N,\widehat h} - b\|_{\delta}^{2} +
\left(\frac{\theta}{2} +\frac{2}{\theta}\right)\|b_{h_0} - b\|_{\delta}^{2} +
\mathfrak c_5\left(\frac{1}{\theta} +\frac{1}{1 -\theta}\right)
\frac{(1 +\lambda)^3}{N}.
\end{displaymath}	
By the decomposition \eqref{risk_bound_PCO_estimator_b_1}, there exist two deterministic constants $\mathfrak c_6,\overline{\mathfrak c}_6 > 0$, not depending on $N$, $\theta$, $h$ and $h_0$, such that with probability larger than $1 -\overline{\mathfrak c}_6|\mathcal H_N|e^{-\lambda}$,
\begin{eqnarray*}
 \|\widehat b_{N,\widehat h} - b\|_{\delta}^{2}
 & \leqslant &
 \|\widehat b_{N,h} - b\|_{\delta}^{2} + |\psi_N(h)| + |\psi_N(\widehat h)|\\
 & \leqslant &
 \left(1 +\frac{\theta}{1 -\theta}\right)
 \|\widehat b_{N,h} - b\|_{\delta}^{2} +
 \frac{\theta}{1 -\theta}\|\widehat b_{N,\widehat h} - b\|_{\delta}^{2}\\
 & &
 \hspace{2cm} +
 \frac{\mathfrak c_6}{\theta}\|b_{h_0} - b\|_{\delta}^{2} +
 \frac{\mathfrak c_6}{\theta(1 -\theta)}\cdot\frac{(1 +\lambda)^3}{N}.
\end{eqnarray*}	
This concludes the proof.
%


%
\subsubsection{Proof of Proposition \ref{kernel_properties_Phi}}\label{Section_proof_kernel_properties_Phi}
For the sake of readability, Proposition \ref{kernel_properties_Phi}.(1,2,3) are proved for $t_0 = 0$ without loss of generality. However, the condition $t_0 > 0$ is required in order to prove Proposition \ref{kernel_properties_Phi}.(4) because $f'$ is involved. First of all, for any $x\in\mathbb R$ and $h > 0$, by It\^o's formula,
\begin{displaymath}
\int_{x_0}^{X_T}K_h(F(z) - F(x))dz =
\int_{0}^{T}K_h(F(X_t) - F(x))dX_t +
\frac{1}{2h^2}\int_{0}^{T}
K'\left(\frac{F(X_t) - F(x)}{h}\right)f(X_t)d\langle X\rangle_t.
\end{displaymath}
So,
\begin{eqnarray*}
 \int_{0}^{T}K_h(F(X_t) - F(x))dX_t & = &
 \int_{x_0}^{X_T}K_h(F(z) - F(x))dz\\
 & & \hspace{1cm} -
 \frac{1}{2h^2}\int_{0}^{T}K'\left(\frac{F(X_t) - F(x)}{h}\right)\sigma(X_t)^2f(X_t)dt
 = T\Phi_h(X,x).
\end{eqnarray*}
\begin{enumerate}
 \item Consider $\varphi\in C^0([0,T];\mathbb R)$. For any $x\in [-\Delta,\Delta]$ and $h\in (0,\Delta_0]$,
 \begin{equation}\label{kernel_properties_Phi_1}
 [-1,1]\subset
 \left[-\frac{F(-\Delta)}{h},\frac{1 - F(\Delta)}{h}\right]\subset
 \left[-\frac{F(x)}{h},\frac{1 - F(x)}{h}\right]
 \end{equation}
 and, since ${\rm supp}(K) = [-1,1]$ and ${\rm supp}(\delta) = [-\Delta,\Delta]$,
 \begin{eqnarray*}
  \int_{-\infty}^{\infty}
  \left[\int_{\varphi(0)}^{\varphi(T)}K_h(F(z) - F(x))dz\right]^2\delta(x)dx
  & \leqslant &
  \int_{-\infty}^{\infty}\left[\int_{-F(x)/h}^{(1 - F(x))/h}
  |K(y)|\frac{dy}{f(F^{-1}(hy + F(x)))}\right]^2\delta(x)dx\\
  & = &
  \int_{-\Delta}^{\Delta}\left(\int_{-1}^{1}
  |K(y)|\frac{dy}{f(F^{-1}(hy + F(x)))}\right)^2\delta(x)dx.
 \end{eqnarray*}
 Moreover, since $h\leqslant\Delta_0$ and $F^{-1}$ is increasing,
 \begin{equation}\label{kernel_properties_Phi_2}
 F^{-1}(h\cdot + F(x))([-1,1])\subset
 [F^{-1}((1 -\kappa)F(-\Delta)),
 F^{-1}(\kappa + (1 -\kappa)F(\Delta))] =: I_{\Delta,F}.
 \end{equation}
 So, there exists $f_1 > 0$ such that $f(z)\geqslant f_1$ for every $z\in I_{\Delta,F}$ because $f(\mathbb R)\subset (0,\infty)$, and then
 \begin{eqnarray*}
  & &
  \int_{-\infty}^{\infty}
  \left[\int_{\varphi(0)}^{\varphi(T)}K_h(F(z) - F(x))dz\right]^2\delta(x)dx\\
  & &
  \hspace{3cm}
  \leqslant
  \frac{1}{f_{1}^{2}}\int_{-\Delta}^{\Delta}
  \left(\int_{-1}^{1}|K(z)|dz\right)^2\delta(x)dx
  \leqslant
  \frac{\|K\|_{1}^{2}}{f_{1}^{2}}.
 \end{eqnarray*}
 The same way, there exists $f_2 > 0$ such that $f(z)\geqslant f_2$ for every $z\in [-\Delta,\Delta]$, and then
 \begin{eqnarray*}
  \int_{-\infty}^{\infty}K'\left(\frac{F(\varphi(t)) - F(x)}{h}\right)^2\delta(x)dx
  & = &
  \int_{-\Delta}^{\Delta}K'\left(\frac{F(\varphi(t)) - F(x)}{h}\right)^2\delta(x)dx\\
  & = &
  h\int_{[F(-\Delta) - F(\varphi(t))]/h}^{[F(\Delta) - F(\varphi(t))]/h}
  K'(y)^2\frac{\delta(F^{-1}(hy + F(\varphi(t))))}{f(F^{-1}(hy + F(\varphi(t))))}dy\\
  & \leqslant &
  \frac{\|K'\|_{2}^{2}\|\delta\|_{\infty}h}{f_2}
 \end{eqnarray*}
 for every $t\in [0,T]$. Thus,
 \begin{eqnarray*}
  T^2\|\Phi_h(\varphi,.)\|_{\delta}^{2}
  & \leqslant &
  2\int_{-\infty}^{\infty}
  \left[\int_{\varphi(0)}^{\varphi(T)}K_h(F(z) - F(x))dz\right]^2\delta(x)dx\\
  & &
  \hspace{2cm} +
  \frac{1}{2h^4}\int_{-\infty}^{\infty}
  \left[\int_{0}^{T}K'\left(\frac{F(\varphi(t)) - F(x)}{h}\right)
  \sigma(\varphi(t))^2f(\varphi(t))dt\right]^2\delta(x)dx\\
  & \leqslant &
  \frac{2\|K\|_{1}^{2}}{f_{1}^{2}} +
  \frac{T}{2h^4}\int_{0}^{T}\sigma(\varphi(t))^4f(\varphi(t))^2
  \int_{-\infty}^{\infty}K'\left(\frac{F(\varphi(t)) - F(x)}{h}\right)^2\delta(x)dxdt\\
  & \leqslant &
  \frac{2\|K\|_{1}^{2}}{f_{1}^{2}} +
  \frac{T^2\|\sigma\|_{\infty}^{4}\|f\|_{\infty}^{2}\|K'\|_{2}^{2}\|\delta\|_{\infty}}{2f_2h^3}.
 \end{eqnarray*}
 \item For any $h,h'\in (0,\Delta_0]$,
 \begin{eqnarray*}
  & &
  \mathbb E(\langle\Phi_h(X^1,.),\Phi_{h'}(X^2,.)\rangle_{\delta}^{2})\\
  & &
  \hspace{1.5cm} =
  \frac{1}{T^4}\mathbb E\left[\left(\int_{-\infty}^{\infty}
  \left(\int_{0}^{T}K_h(F(X_{t}^{1}) - F(x))dX_{t}^{1}\right)
  \left(\int_{0}^{T}K_{h'}(F(X_{t}^{2}) - F(x))dX_{t}^{2}\right)\delta(x)dx\right)^2
  \right]\\
  & &
  \hspace{1.5cm}\leqslant
  \frac{2}{T^4}(\mathbb E(A_{h,h'}^{2})
  +\mathbb E(B_{h,h'}^{2}))
 \end{eqnarray*}
 with
 \begin{displaymath}
 A_{h,h'} :=
 \int_{-\infty}^{\infty}
 \left(\int_{0}^{T}K_h(F(X_{t}^{1}) - F(x))\sigma(X_{t}^{1})dW_{t}^{1}\right)
 \left(\int_{0}^{T}K_{h'}(F(X_{t}^{2}) - F(x))dX_{t}^{2}\right)\delta(x)dx
 \end{displaymath}
 and
 \begin{displaymath}
 B_{h,h'} :=
 \int_{-\infty}^{\infty}
 \left(\int_{0}^{T}K_h(F(X_{t}^{1}) - F(x))b(X_{t}^{1})dt\right)
 \left(\int_{0}^{T}K_{h'}(F(X_{t}^{2}) - F(x))dX_{t}^{2}\right)
 \delta(x)dx.
 \end{displaymath}
 {\bf Bound on $\mathbb E(A_{h,h'}^{2})$.} Since $(X^1,W^1)$ and $X^2$ are independent,
 \begin{eqnarray*}
  \mathbb E(A_{h,h'}^{2}) & = &
  \int_{-\infty}^{\infty}\int_{-\infty}^{\infty}
  \mathbb E\left[\left(\int_{0}^{T}K_h(F(X_{t}^{1}) - F(x))\sigma(X_{t}^{1})dW_{t}^{1}\right)
  \left(\int_{0}^{T}K_h(F(X_{t}^{1}) - F(y))\sigma(X_{t}^{1})dW_{t}^{1}\right)\right]\\
  & &
  \hspace{1cm}\times
  \mathbb E\left[\left(\int_{0}^{T}K_{h'}(F(X_{t}^{2}) - F(x))dX_{t}^{2}\right)
  \left(\int_{0}^{T}K_{h'}(F(X_{t}^{2}) - F(y))dX_{t}^{2}\right)\right]\delta(x)\delta(y)dxdy.
 \end{eqnarray*}
 On the one hand, for every $x,y\in\mathbb R$, by the isometry property of It\^o's integral and the definition of $f$,
 \begin{eqnarray*}
  & &
  \mathbb E\left[\left(\int_{0}^{T}K_h(F(X_{t}^{1}) - F(x))\sigma(X_{t}^{1})dW_{t}^{1}\right)
  \left(\int_{0}^{T}K_h(F(X_{t}^{1}) - F(y))\sigma(X_{t}^{1})dW_{t}^{1}\right)\right]\\
  & &
  \hspace{3cm}
  =\int_{0}^{T}\mathbb E(K_h(F(X_{t}^{1}) - F(x))
  K_h(F(X_{t}^{1}) - F(y))\sigma(X_{t}^{1})^2)dt\\
  & &
  \hspace{6cm}
  = T\int_{-\infty}^{\infty}K_h(F(z) - F(x))K_h(F(z) - F(y))\sigma(z)^2f(z)dz.
 \end{eqnarray*}
 Then,
 \begin{eqnarray*}
  \mathbb E(A_{h,h'}^{2}) & = &
  T\int_{-\infty}^{\infty}\int_{-\infty}^{\infty}
  \int_{-\infty}^{\infty}K_h(F(z) - F(x))K_h(F(z) - F(y))\sigma(z)^2f(z)\\
  & &
  \hspace{0.5cm}\times
  \mathbb E\left[
  \left(\int_{0}^{T}K_{h'}(F(X_{t}^{2}) - F(x))dX_{t}^{2}\right)
  \left(\int_{0}^{T}K_{h'}(F(X_{t}^{2}) - F(y))dX_{t}^{2}\right)\right]\delta(x)\delta(y)dxdydz\\
  & = &
  T\int_{-\infty}^{\infty}\sigma(z)^2f(z)
  \mathbb E\left[\left(\int_{-\Delta}^{\Delta}
  K_h(F(z) - F(x))\delta(x)
  \int_{0}^{T}K_{h'}(F(X_{t}^{2}) - F(x))dX_{t}^{2}dx\right)^2\right]dz.
 \end{eqnarray*}
 On the other hand, note that for any $z\in\mathbb R$,
 \begin{eqnarray}
  \|K_h(F(z) - F(.))\mathbf 1_{[-\Delta,\Delta]}(.)\|_1 & = &
  \int_{-\Delta}^{\Delta}|K_h(F(z) - F(y))|dy
  \nonumber\\
  \label{kernel_properties_Phi_3}
  & = &
  \int_{(F(-\Delta) - F(z))/h}^{(F(\Delta) - F(z))/h}
  |K(y)|\frac{dy}{f(F^{-1}(hy + F(z)))}\leqslant\frac{\|K\|_1}{f_2}
 \end{eqnarray}
 and, for every $x\in\mathbb R$,
 \begin{displaymath}
 \int_{-\infty}^{\infty}f(z)|K_h(F(z) - F(x))|dz =
 \int_{0}^{1}|K_h(z - F(x))|dz\leqslant\|K\|_1.
 \end{displaymath}
 Then, since
 \begin{displaymath}
 x\in [-\Delta,\Delta]\longmapsto
 \frac{|K_h(F(z) - F(x))|}{\|K_h(F(z) - F(.))\mathbf 1_{[-\Delta,\Delta]}(.)\|_1}
 \end{displaymath}
 is a density function, by Jensen's inequality,
 \begin{eqnarray*}
  \mathbb E(A_{h,h'}^{2}) & \leqslant &
  T\int_{-\infty}^{\infty}\sigma(z)^2f(z)\|K_h(F(z) - F(.))\mathbf 1_{[-\Delta,\Delta]}(.)\|_1\\
  & &
  \hspace{1cm}
  \times\int_{-\Delta}^{\Delta}
  |K_h(F(z) - F(x))|\delta(x)^2\mathbb E\left[\left(
  \int_{0}^{T}K_{h'}(F(X_{t}^{2}) - F(x))dX_{t}^{2}\right)^2\right]dxdz\\
  & \leqslant &
  \frac{T\|\sigma\|_{\infty}^{2}\|K\|_{1}^{2}\|\delta\|_{\infty}}{f_2}
  \int_{-\Delta}^{\Delta}\delta(x)\mathbb E\left[\left(
  \int_{0}^{T}K_{h'}(F(X_{t}^{2}) - F(x))dX_{t}^{2}\right)^2\right]dx\\
  & &
  \hspace{8cm}\leqslant
  \frac{T^3\|\sigma\|_{\infty}^{2}
  \|K\|_{1}^{2}\|\delta\|_{\infty}\mathfrak m(h')}{f_2}.
 \end{eqnarray*}
 {\bf Bound on $\mathbb E(B_{h,h'}^{2})$.} First, note that for every $x\in [-\Delta,\Delta]$, by (\ref{kernel_properties_Phi_1}) and (\ref{kernel_properties_Phi_2}),
 \begin{eqnarray}
  \int_{-\infty}^{\infty}|K_h(F(z) - F(x))|b(z)^2f(z)dz & = &
  \int_{-F(x)/h}^{(1 - F(x))/h}|K(z)|b(F^{-1}(hz + F(x)))^2dz
  \nonumber\\
  \label{kernel_properties_Phi_4}
  & = &
  \int_{-1}^{1}|K(z)|b(F^{-1}(hz + F(x)))^2dz\leqslant
  \|b\|_{\infty,I_{\Delta,F}}^{2}\|K\|_1.
 \end{eqnarray}
 Then, since
 \begin{displaymath}
 x\in [-\Delta,\Delta]\longmapsto
 \frac{|K_h(F(X_t(\omega)) - F(x))|}{\|K_h(F(X_t(\omega)) - F(.))
 \mathbf 1_{[-\Delta,\Delta]}(.)\|_1}
 \end{displaymath}
 is a density function for every $(t,\omega)\in [0,T]\times\Omega$, by Jensen's inequality and Inequality (\ref{kernel_properties_Phi_3}),
 \begin{eqnarray*}
  \mathbb E(B_{h,h'}^{2}) & = &
  \mathbb E\left[\left(\int_{0}^{T}\int_{-\infty}^{\infty}K_h(F(X_{t}^{1}) - F(x))
  b(X_{t}^{1})\delta(x)
  \int_{0}^{T}K_{h'}(F(X_{s}^{2}) - F(x))dX_{s}^{2}dxdt\right)^2\right]\\
  & \leqslant &
  \frac{T\|K\|_1}{f_2}
  \int_{0}^{T}\int_{-\infty}^{\infty}\mathbb E(|K_h(F(X_{t}^{1}) - F(x))|b(X_{t}^{1})^2)\\
  & &
  \hspace{4cm}\times
  \delta(x)^2\mathbb E\left[
  \left(\int_{0}^{T}K_{h'}(F(X_{s}^{2}) - F(x))dX_{s}^{2}\right)^2\right]dxdt\\
  & = &
  \frac{T^2\|K\|_1}{f_2}
  \int_{-\infty}^{\infty}\left(\int_{-\infty}^{\infty}
  |K_h(F(z) - F(x))|b(z)^2f(z)dz\right)\\
  & &
  \hspace{4cm}\times
  \delta(x)^2\mathbb E\left[
  \left(\int_{0}^{T}K_{h'}(F(X_{s}^{2}) - F(x))dX_{s}^{2}\right)^2\right]dx\\
  & \leqslant &
  \frac{T^2\|b\|_{\infty,I_{\Delta,F}}^{2}\|K\|_{1}^{2}\|\delta\|_{\infty}}{f_2}
  \int_{-\infty}^{\infty}\delta(x)\mathbb E\left[
  \left(\int_{0}^{T}K_{h'}(F(X_{s}^{2}) - F(x))dX_{s}^{2}\right)^2\right]dx\\
  & &
  \hspace{8cm}\leqslant
  \frac{T^4\|b\|_{\infty,I_{\Delta,F}}^{2}\|K\|_{1}^{2}\|\delta\|_{\infty}\mathfrak m(h')}{f_2}.
 \end{eqnarray*}
 \item For any $h\in (0,\Delta_0]$ and $\varphi\in\mathbb L^2(\mathbb R,dx)$,
 \begin{eqnarray*}
  \mathbb E(\langle\Phi_h(X,.),\varphi\rangle_{\delta}^{2}) & = &
  \frac{1}{T^2}
  \mathbb E\left[\left(\int_{-\infty}^{\infty}
  \varphi(x)\delta(x)\int_{0}^{T}K_h(F(X_t) - F(x))dX_tdx\right)^2
  \right]\\
  & \leqslant &
  \frac{2}{T^2}(\mathbb E(C_{h}^{2}) +\mathbb E(D_{h}^{2}))
 \end{eqnarray*}
 with
 \begin{displaymath}
 C_h :=
 \int_{-\infty}^{\infty}
 \varphi(x)\delta(x)\int_{0}^{T}K_h(F(X_t) - F(x))\sigma(X_t)dW_tdx
 \end{displaymath}
 and
 \begin{displaymath}
 D_h :=
 \int_{-\infty}^{\infty}
 \varphi(x)\delta(x)\int_{0}^{T}K_h(F(X_t) - F(x))b(X_t)dtdx.
 \end{displaymath}
 {\bf Bound on $\mathbb E(C_{h}^{2})$.} By the isometry property of It\^o's integral and Jensen's inequality,
 \begin{eqnarray*}
  \mathbb E(C_{h}^{2}) & = &
  \int_{-\infty}^{\infty}\int_{-\infty}^{\infty}
  \varphi(x)\varphi(y)\delta(x)\delta(y)\int_{0}^{T}
  \mathbb E(K_h(F(X_t) - F(x))K_h(F(X_t) - F(y))
  \sigma(X_t)^2)dtdxdy\\
  & = &
  T\int_{-\infty}^{\infty}\int_{-\infty}^{\infty}
  \varphi(x)\varphi(y)\delta(x)\delta(y)\\
  & &
  \hspace{3cm}\times
  \int_{-\infty}^{\infty}K_h(F(z) - F(x))K_h(F(z) - F(y))\sigma(z)^2f(z)dzdxdy\\
  & = &
  T\int_{-\infty}^{\infty}\left(\int_{-\Delta}^{\Delta}
  K_h(F(z) - F(x))\varphi(x)\delta(x)dx\right)^2\sigma(z)^2f(z)dz\\
  & \leqslant &
  T\int_{-\infty}^{\infty}\|K_h(F(z) - F(.))\mathbf 1_{[-\Delta,\Delta]}(.)\|_1\\
  & &
  \hspace{3cm}\times
  \left(\int_{-\Delta}^{\Delta}
  |K_h(F(z) - F(x))|\varphi(x)^2\delta(x)^2dx\right)\sigma(z)^2f(z)dz\\
  & \leqslant &
  \frac{T\|K\|_1\|\sigma\|_{\infty}^{2}}{f_2}\int_{-\Delta}^{\Delta}\left(
  \int_{-\infty}^{\infty}f(z)|K_h(F(z) - F(x))|dz\right)
  \varphi(x)^2\delta(x)^2dx
  \quad\textrm{by Inequality (\ref{kernel_properties_Phi_3})}\\
  & \leqslant &
  \frac{T\|K\|_{1}^{2}\|\sigma\|_{\infty}^{2}
  \|\delta\|_{\infty}\|\varphi\|_{\delta}^{2}}{f_2}.
 \end{eqnarray*}
 {\bf Bound on $\mathbb E(D_{h}^{2})$.} By Jensen's inequality, and Inequalities (\ref{kernel_properties_Phi_3}) and (\ref{kernel_properties_Phi_4}),
 \begin{eqnarray*}
  \mathbb E(D_{h}^{2})
  & \leqslant &
  T\int_{0}^{T}\mathbb E\left[b(X_t)^2\|K_h(F(X_t) - F(.))\mathbf 1_{[-\Delta,\Delta]}(.)\|_1
  \int_{-\Delta}^{\Delta}|K_h(F(X_t) - F(x))|\varphi(x)^2\delta(x)^2dx\right]dt\\
  & \leqslant &
  \frac{T^2\|K\|_1}{f_2}\int_{-\Delta}^{\Delta}
  \left(\int_{-\infty}^{\infty}|K_h(F(z) - F(x))|b(z)^2f(z)dz\right)\varphi(x)^2\delta(x)^2dx\\
  & &
  \hspace{8cm}\leqslant
  \frac{T^2\|b\|_{\infty,I_{\Delta,F}}^{2}\|K\|_{1}^{2}
  \|\delta\|_{\infty}\|\varphi\|_{\delta}^{2}}{f_2}.
 \end{eqnarray*}
 \item Since $X$ is a semi-martingale, since the map
 \begin{displaymath}
 (t,\omega,x)\longmapsto
 K_h(F(X_t(\omega)) - F(x))b_{h'}(x)\mathbf 1_{[-\Delta,\Delta]}(x)
 \end{displaymath}
 is measurable and bounded for any $h,h'\in\mathcal H_N$, and since
 \begin{displaymath}
 A\mapsto\int_A\delta(x)dx
 \quad\textrm{is a finite measure,}
 \end{displaymath}
 by the stochastic Fubini theorem, the change of variable formula and It\^o's formula,
 \begin{eqnarray*}
  (T - t_0)\langle\Phi_h(X,.),b_{h'}\rangle_{\delta}
  & = &
  (T - t_0)\int_{-\infty}^{\infty}
  \Phi_h(X,x)b_{h'}(x)\delta(x)dx\\
  & = &
  \int_{t_0}^{T}
  \int_{-\infty}^{\infty}
  K_h(F(X_t) - F(x))b_{h'}(x)\delta(x)dxdX_t
  \quad {\rm a.s.}\\
  & = &
  \int_{t_0}^{T}[K_h\ast [((b_{h'}\delta/f)\circ F^{-1})\mathbf 1_{(0,1)}]](F(X_t))dX_t\\
  & = &
  \int_{X_{t_0}}^{X_T}\psi_{h,h'}(z)dz -
  \frac{1}{2}\int_{t_0}^{T}\psi_{h,h'}'(X_t)\sigma(X_t)^2dt
 \end{eqnarray*}
 where
 \begin{displaymath}
 \psi_{h,h'}(.) := [K_h\ast [((b_{h'}\delta/f)\circ F^{-1})\mathbf 1_{(0,1)}]](F(.)).
 \end{displaymath}
 On the one hand,
 \begin{eqnarray*}
  (\psi_{h,h'}'/f)\circ F^{-1} & = &
  K_h\ast [((b_{h'}\delta/f)\circ F^{-1})\mathbf 1_{(0,1)}]'\\
  & = &
  K_h\ast [\mathbf 1_{[F(-\Delta),F(\Delta)]}/(f\circ F^{-1})\times
  (b_{h'}\delta'/f + b_{h'}'\delta/f - b_{h'}\delta f'/f^2)\circ F^{-1}].
 \end{eqnarray*}
 Moreover, for every $x\in\mathbb R$, by (\ref{kernel_properties_Phi_2}),
 \begin{eqnarray}
  |b_{h'}(x)\mathbf 1_{[-\Delta,\Delta]}(x)|
  & = &
  \left|\mathbf 1_{[-\Delta,\Delta]}(x)
  \int_{-\infty}^{\infty}K_{h'}(y - F(x))[(b\circ F^{-1})\mathbf 1_{(0,1)}](y)dy\right|
  \nonumber\\
  \label{kernel_properties_Phi_5}
  & \leqslant &
  \|K\|_1\sup_{(x,y)\in [-\Delta,\Delta]\times [-1,1]}
  |[(b\circ F^{-1})\mathbf 1_{(0,1)}](h'y + F(x))|
  \leqslant
  \|K\|_1\|b\|_{\infty,I_{\Delta,F}}
 \end{eqnarray}
 and
 \begin{eqnarray*}
  |b_{h'}'(x)\mathbf 1_{[-\Delta,\Delta]}(x)|
  & = &
  \left|\mathbf 1_{[-\Delta,\Delta]}(x)f(x)
  \int_{-\infty}^{\infty}K_{h'}(y - F(x))[(b\circ F^{-1})'\mathbf 1_{(0,1)}](y)dy\right|\\
  & \leqslant &
  \|K\|_1\sup_{(x,y)\in [-\Delta,\Delta]\times [-1,1]}
  \left|\left[\frac{b'\circ F^{-1}}{f\circ F^{-1}}\mathbf 1_{(0,1)}\right](h'y + F(x))\right|
  \leqslant
  \frac{\|K\|_1\|b'\|_{\infty,I_{\Delta,F}}}{f_1}.
 \end{eqnarray*}
 Then, since $f$ and $f'$ are bounded by Remark \ref{F_properties},
 \begin{eqnarray*}
  \|\psi_{h,h'}'\|_{\infty}
  & \leqslant &
  \frac{\|f\|_{\infty}\|K_h\|_1}{f_2}\left(
  \frac{\|\delta'\|_{\infty}\|b_{h'}\mathbf 1_{[-\Delta,\Delta]}\|_{\infty}}{f_2} +
  \frac{\|\delta\|_{\infty}\|b_{h'}'\mathbf 1_{[-\Delta,\Delta]}\|_{\infty}}{f_2} +
  \frac{\|\delta\|_{\infty}\|f'\|_{\infty}\|b_{h'}\mathbf 1_{[-\Delta,\Delta]}\|_{\infty}}{f_{2}^{2}}
  \right)\\
  & \leqslant &
  \frac{\|f\|_{\infty}\|K\|_{1}^{2}}{f_2}\left(
  \frac{\|\delta'\|_{\infty}\|b\|_{\infty,I_{\Delta,F}}}{f_2} +
  \frac{\|\delta\|_{\infty}\|b'\|_{\infty,I_{\Delta,F}}}{f_1f_2} +
  \frac{\|\delta\|_{\infty}\|f'\|_{\infty}\|b\|_{\infty,I_{\Delta,F}}}{f_{2}^{2}}
  \right) <\infty.
 \end{eqnarray*}
 On the other hand, by Inequality (\ref{kernel_properties_Phi_5}), by (\ref{kernel_properties_Phi_1}) and by (\ref{kernel_properties_Phi_2}),
 \begin{eqnarray*}
  \int_{-\infty}^{\infty}|\psi_{h,h'}(x)|dx
  & = &
  \int_{-\infty}^{\infty}\left|\int_{0}^{1}K_h(F(x) - y)
  \left(\frac{b_{h'}\delta}{f}\right)(F^{-1}(y))dy\right|dx\\
  & = &
  \int_{0}^{1}\left|\int_{-\infty}^{\infty}K_h(x - F(y))
  b_{h'}(y)\delta(y)dy\right|\frac{dx}{f(F^{-1}(x))}\\
  & \leqslant &
  \|K\|_1\|b\|_{\infty,I_{\Delta,F}}
  \int_{0}^{1}\int_{-\Delta}^{\Delta}|K_h(x - F(y))|
  \delta(y)\frac{dydx}{f(F^{-1}(x))}\\
  & = &
  \|K\|_1\|b\|_{\infty,I_{\Delta,F}}
  \int_{-\Delta}^{\Delta}\delta(y)
  \int_{-F(y)/h}^{(1 - F(y))/h}|K(x)|
  \frac{dx}{f(F^{-1}(hx + F(y)))}dy\\
  & \leqslant &
  \frac{\|K\|_1\|b\|_{\infty,I_{\Delta,F}}}{f_1}
  \int_{-\Delta}^{\Delta}\delta(y)
  \int_{-1}^{1}|K(x)|dxdy =
  \frac{\|K\|_{1}^{2}\|b\|_{\infty,I_{\Delta,F}}}{f_1}.
 \end{eqnarray*}
 This concludes the proof because
 \begin{displaymath}
 (T - t_0)|\langle\Phi_h(X,.),b_{h'}\rangle_{\delta}|
 \leqslant
 \int_{-\infty}^{\infty}|\psi_{h,h'}(x)|dx +
 \frac{T - t_0}{2}\|\sigma\|_{\infty}^{2}\|\psi_{h,h'}'\|_{\infty}
 \quad {\rm a.s.}
 \end{displaymath}
\end{enumerate}
%


%

%
\end{document}